\documentclass[10pt]{amsart}

\usepackage{epic}
\usepackage{eepic}
\usepackage{graphicx, color}
\usepackage{psfrag}

\graphicspath{{figures/}}

\newtheorem{thm}{Theorem}[section]

\newtheorem{coro}[thm]{Corollary}
\newtheorem{lem}[thm]{Lemma}
\newtheorem{prop}[thm]{Proposition}

\theoremstyle{definition}

\theoremstyle{remark}

\newtheorem{remk}[thm]{Remark}

\newcommand{\R}{{\mathbb{R}}}

\newcommand{\refeq}[1]{(\ref{#1})}

\begin{document}

\title[Generic spectral simplicity of polygons]{Generic spectral simplicity of polygons}

\author{L. Hillairet and C. Judge}     

\thanks{We thank the Fields Institute for hosting 
a workshop where this work began.   
We also thank Gilles  Carron and David Hoff for
helpful discussions.}

\date{\today}

\maketitle

\begin{abstract}
We study the Laplace operator with Dirichlet or Neumann 
boundary condition on polygons in the Euclidean 
plane. We prove that almost every simply connected polygon with at least four 
vertices has simple spectrum. We also address the more
general case of geodesic polygons in a constant curvature space form.
%
%
\end{abstract}
  
\section{Introduction}

Let ${\mathcal P}_n$ be the space of 
polygons in ${\mathbb R}^2$ having $n$-vertices and  
boundary equal to a simple closed curve. 
By labeling the vertices, one obtains a local 
parametrization of ${\mathcal P}_n$ by ${\mathbb R}^{2n}$. 
This defines an affine structure on ${\mathcal P}_n$
as well as a notion of sets of measure zero.

\begin{thm} \label{TheoremSimple}
If $n \geq 4$, then for almost every $P \in {\mathcal P}_n$, 
the Dirichlet (resp. Neumann) Euclidean Laplacian on $P$ 
has one-dimensional eigenspaces. 
\end{thm} 

As of yet, we do not know whether or not
the generic triangle has simple spectrum.  
The methods of this paper show that either
the generic triangle has simple spectrum or 
no triangle has simple spectrum.

We also generalize Theorem \ref{TheoremSimple} to 
\begin{thm}\label{TheoremConstCurvPoly}
Almost every simply connected geodesic polygon 
belonging to a hemisphere of the standard sphere  has simple 
Dirichlet (resp. Neumann)  spectrum.
Almost every simply connected geodesic polygon in the 
hyperbolic plane has simple Dirichlet (resp. Neumann) 
spectrum.  
\end{thm} 
It is conjectured by number theorists that the geodesic 
triangle in the hyperbolic plane ${\mathbb H}^2$
with angles $\pi/2$, $\pi/6$, and $0$, has simple spectrum 
\cite{Srn}. The (joint) spectrum of this triangle 
coincides with the spectrum of the hyperbolic surface
${\mathbb H}^{2} / {\rm SL}(2, {\mathbb Z})$. 
It is not known whether there exists a closed
hyperbolic surface with simple Laplace spectrum.

We note that K. Uhlenbenbeck \cite{Uhlenbeck}---building 
upon earlier work of J. Albert \cite{A}---proved the simplicity
of the Laplace spectrum for the generic compact domain in 
${\mathbb R}^d$ and the generic metric on a compact
manifold.  Such moduli spaces are infinite dimensional 
whereas ${\mathcal P}_n$ is only finite dimensional. 


\section{Analytic bi-Lipschitz boundary perturbations} 

The proofs of Theorems  \ref{TheoremSimple} and 
\ref{TheoremConstCurvPoly} are based upon 
analytic perturbation theory. 
We first show that the Laplace spectrum of a real-analytic 
path $P_t$ of polygons depends real-analytically on $t$. 
This will follow from   

\begin{lem}  \label{AnalyticEigenvalue} 
Let $\Omega_t$, $ 0 \leq t \leq 1$,  
be a family of Lipschitz domains such that there
exists a real-analytic family of homeomorphisms 
$f_t: \Omega_0 \rightarrow \Omega_t$ and a constant $C>0$
such that 
\begin{equation} \label{LipschitzBound}
   C\cdot 
 |x-y|~ \leq~ |f_t(x) - f_t(y)|~ \leq~  \frac{1}{C} \cdot |x-y|. 
\end{equation}
Then the eigenvalues of Neumann (resp. Dirichlet) 
Laplacian, $\Delta_t$, vary real-analytically in $t$. 
\end{lem}
 
\begin{proof}\footnote{
If $f_t$ is a real-analytic family of smooth 
diffeomorphisms, this fact belongs to the 
standard theory of {\it boundary perturbation}. 
See, for instance, Kato \S VII.6.5. 
Our proof will follow the same lines.}
For each $t_0$, there exists an interval neighborhood 
$I \ni t_0$ and Lipschitz functions $f^{(k)}$ such that 
$f_t(x) =  \sum_k  f^{(k)}(x) (t-t_0)^k$. 
By Rademacher's theorem, for each $k$, there exists a 
set of full measure $A_k \subset \Omega_0$ such that 
for each $x \in \Omega_0$, the differential $D_x f^{(k)}$ 
exists. Letting $A = \cap_k  A_k$ we obtain a set of full 
measure such that $D_x f_t$ exists for all $x \in A$ and 
$t \in I$.

From (\ref{LipschitzBound})  we have 
$C \leq |D_x f_t| \leq C^{-1}$ for $x \in A$.  
It follows that the pull-back operator 
$f_t^*: C^{\infty}(\Omega_t) \rightarrow C^{\infty}(\Omega_0)$
extends continuously in the $H^1$ norm to an 
isomorphism from  $H^1(\Omega_t)$ to $H^1(\Omega_0)$. 
Also, it follows that the measure
$(f_t)^{-1}_*(\mu)$ is absolutely continuous with respect
to Lebesgue measure, and hence equals $h_t d \mu$ for 
some a.e. positive function 
$h_t: \Omega_t \rightarrow {\mathbb R}$. 

For any $v \in H^1(\Omega_t),$ we then have that 
\begin{gather}\label{CV1}
\int_{\Omega_t} |v(X)|^2 d\mu(X)\,=\,\int_{\Omega_0} |v\circ f_t(x)|^2 h_t(x) d\mu(x)\\
\label{CV2} \nabla_x (v\circ f_t) \,=\, D_xf_t \cdot \left[ (\nabla_X v)\circ f_t\right ] .
\end{gather}
It follows that the operator $H_t$ 
defined by 
\[     H_t(u)~  =~ (u \circ f_t) \]
is unitary  when acting from $L^2(\Omega_t, d\mu)$ 
onto $L^2(\Omega_0, h_t d\mu).$ 

Let us consider first the Neumann eigenvalue problem i.e. $\Delta_t$ is the 
self-adjoint operator associated with the quadratic form $q_t$ defined on $H^1(\Omega_t).$
Define the quadratic form $\tilde{q}_t$ on $H^1(\Omega_0)$ 
by $\tilde{q}_t =  q_t \circ H_t^{-1}.$ 
Using \refeq{CV1} and \refeq{CV2} we find that
$$ {\tilde q}_t(u)\,=\,\int_{\Omega_0} | (D_xf_t)^{-1} \nabla_x u|^2 \,h_t d\mu.$$
Denote by  $\tilde{\Delta}_t$ the associated self-adjoint operator 
(with respect to $L^2(\Omega_0, h_t d\mu)$).  
By definition, $H_t$ defines a unitary equivalence, and in particular  
$\tilde{\Delta}_t$ and  $\Delta_t$ have the same Neumann spectrum.  

Using the power series expansion for $f_t$ in $t$, we find 
that, for each $u \in V_0$, the maps $t \mapsto {\tilde q}_t(u)$ and 
$t \mapsto \int_{\Omega_0} |u|^2\, h_t d\mu$ are real-analytic. 
The spectrum of $\tilde{\Delta}_t$ is thus determined by a generalized 
eigenvalue problem in the sense of Kato (sect. VII.6.4-5). 
The conclusion thus holds.

In the case of Dirichlet eigenvalues, 
$\Delta_t$ is the self-adjoint operator 
associated with $q_t$ defined on $H_0^1(\Omega_t)$,  the 
closure in $H^1(\Omega_t)$ of the smooth functions on $\Omega_t$
that vanish on the boundary of $\Omega_t$. Since $H_t$ maps 
$H^1_0(\Omega_t)$ onto $H^1_0(\Omega_0)$, the same method
of proof applies.
\end{proof}

\begin{coro} \label{LemmaAnalyticSimple}
If for some $t_0$, the spectrum of $\Delta_{t_0}$
is simple, then for all but countably many $t$,
the spectrum of  $\Delta_{t}$ is simple. 
\end{coro}

\begin{proof}
Let $\lambda_n(t)$ be the real-analytic functions
corresponding to the eigenvalues of $\Delta_t$.
Since the spectrum of $\Delta_{t_0}$ is simple,
$\lambda_i(t_0) \neq \lambda_j(t_0)$ if $ i \neq j$.
Thus, since they are real-analytic,  the set 
$\{ t ~|~ \lambda_i(t) = \lambda_j(t) \}$ is at most
countable for every $i \neq j$.  A countable union
of countable sets is countable. 
\end{proof}

In the sequel, we will consider
linear families of piecewise linear homeomorphisms.
These are real-analytic families of
uniformly bi-Lipschitz homeomorphisms that describe 
all deformations of polygons.  
Corollary 2.2 thus implies the following proposition.

\begin{prop} \label{SimpleNeighborhood}
If the Laplace spectrum of $P \in {\mathcal P}^n$ 
is simple, then there exists an open ball neighborhood
${\mathcal B}$ of 
$P \in {\mathcal P}_n \subset {\mathbb R}^{2n}$ 
such that almost every polygon in $Q \in {\mathcal B}$
has simple spectrum. 
\end{prop}

\begin{proof} 
Triangulate $P$ so that the vertices of the triangulation 
that lie in $\partial P$ are exactly the vertices of $P$.
Let $v_{n+1}, \cdots, v_{n+k}$ denote the additional vertices.
Let ${\mathcal B} \subset {\mathbb R}^{2n}$ be a ball 
neighborhood of $P= (v_1, v_2, \ldots, v_n)$ such that 
for each $P' = (v_1', v_2', \ldots, v_n') \in {\mathcal B}$,
we have that $P'$ is a simply connected polygon and 
$v_1', v_2', \cdots, v_n'$,   $ v_{n+1}, \cdots v_{n+k}$ 
are the vertices of a triangulation of $P'$. 
Using the triangulation, we can define a piecewise 
linear homeomorphism $f$ such that $P'=f(P)$.   
By letting $f_t= (1-t){\rm Id} +t f$ and $P_t=f_t(P)$, 
we obtain a linear family of piecewise linear homeomorphisms
corresponding to the segment in ${\mathcal P}^n$
that joins $P$ to $P'$.

Thus, by Corollary \ref{LemmaAnalyticSimple},
for almost all $P_t$ in the segment joining $P$ to $P'$, 
the spectrum of $\Delta_Q$ is simple. The claim 
then follows from using polar coordinates based at $P$
and applying Fubini's theorem.
\end{proof}


\section{Proof of Theorem \ref{TheoremSimple}}

Proposition \ref{SimpleNeighborhood} reduces the proof
of Theorem  \ref{TheoremSimple} to verifying that for each 
$n \geq 4$, the set ${\mathcal P}^n \subset {\mathbb R}^{2n}$
is connected  and contains at least one polygon $P$ with
simple spectrum. We first show that every for $n \geq 4$
there exists an $n$-gon whose Laplace spectrum is simple.

Let $R$ be a rectangle with sidelengths $s_1$ and $s_2$.
An $L^2(R)$ basis of Neumann (resp. Dirichlet) Laplace 
eigenfunctions of $R$
is given by $\cos( \pi m x/s_1) \cdot \cos(\pi n x/s_2)$
(resp. $\sin( \pi m x/s_1) \cdot \sin(\pi n x/s_2)$) where 
$m$ and $n$ vary over the nonnegative integers 
(resp. positive integers).  In particular, we have
multiple eigenvalues iff there exist pairs of integers
$(m, n)$  and $(\bar{m},\bar{n})$ so that    

\[  \frac{m^2}{s_1^2} + \frac{n^2}{s_2^2}~ =~
    \frac{\bar{m}^2}{s_1^2} + \frac{\bar{n}^2}{s_2^2}.  \]
Thus if $(s_1/s_2)^2 \notin {\mathbb Q}$, then the
spectrum of the rectangle $R$ is simple. 

By adding $n-4$ vertices to a side of a rectangle,
one may regard any rectangle as an $n$-gon. 
The spectrum is unchanged by adding such `false' vertices.
Thus, for any $n \geq 4$, there exists an $n$-gon
with simple spectrum. 

The (path) connectedness of ${\mathcal P}^n$ can be 
verified by induction on the number of vertices \footnote{We remark 
that any star-shaped 
$n$-gon can be easily connected to a rectangle with $n-4$ `false' vertices.} 
The case of 3 vertices
can be verified in many ways.  For example,  
the connected Lie group consisting of affine 
homeomorphisms of ${\mathbb R}^2$
acts transitively and continously
on triples of distinct points in ${\mathbb R}^2$.  

Recall that an $n+1$-gon may be regarded as an $n$-gon
if three successive vertices lie on the same line
segment. We will show that any $n+1$-gon, $n \geq 3$, can 
be linearly deformed into an $n$-gon. Thus, the 
claim will follow from the induction hypothesis.

Any $n+1$-gon $P=(v_1, \ldots, v_{n+1})$
can be triangulated so that the vertices of 
the triangulation are exactly the vertices
$v_1, \ldots, v_{n+1}$ of the polygon.\footnote{For example, 
the Delaunay triangulation of $P$ has this property \cite{Thr}.}
Since $P$ is simply connected and $n+1 \geq 3$, 
the dual graph of this triangulation is a tree with 
at least two vertices. Let $T$ be a triangle 
corresponding to an end of the dual graph.
Then there is a vertex $v$ of $T$ such that
both of the sides of $T$ adjacent to $v$
are also sides of $P$. Without loss of generality
$v=v_1$. 

Let $m$ be the midpoint (for example)
of the side of $T$ that is opposite to $v_i$.
Define 
\[    P_t~  =~  \left( (1-t) \cdot v_1 + t \cdot m,
     v_2, \ldots, v_{n+1} \right). \]
Note that the vertex $m$ lies in a segment 
joining two other vertices of $P_1= (m, v_2,\ldots, v_n)$. 
Hence $P_1$ can be regarded as an $n$-gon.

\begin{figure}[h]
\begin{center}

  \psfrag{v1}{$v_1$}
  \psfrag{v2}{$v_2$}
  \psfrag{vn+1}{$v_{n+1}$}
  \psfrag{m}{$m$}

\includegraphics{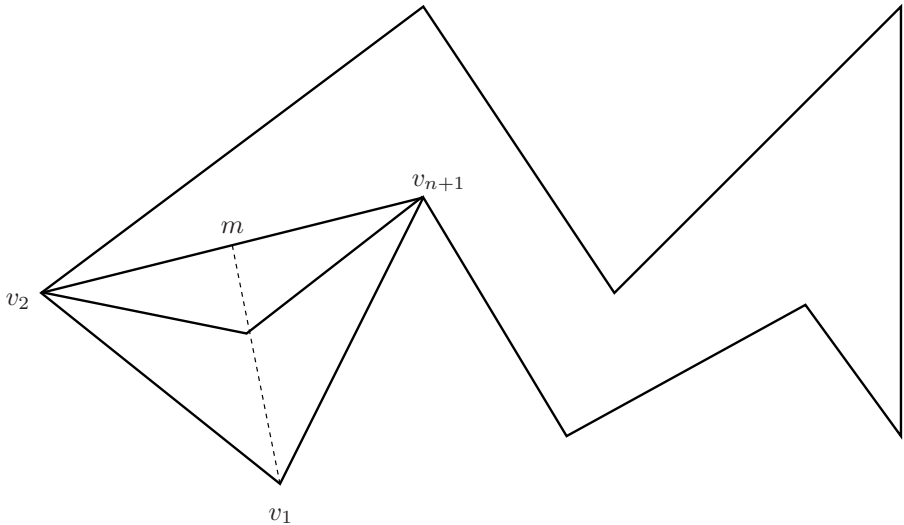}%

\end{center}
\caption{\label{FigDeleting}Deleting a vertex}
\end{figure}


\section{Proof of Theorem \ref{TheoremConstCurvPoly}}

In this section we prove Theorem \ref{TheoremConstCurvPoly}. 
The key is to construct, for each Euclidean polygon $P_0$, 
a natural real-analytic family of polygons 
$\kappa \rightarrow P_{\kappa}$ such that $P_{\kappa}$ 
is a geodesic polygon in the constant curvature $\kappa$ 
space form
$M_{\kappa}$.  To make the construction transparent, we 
make a convenient choice of model for $M_{\kappa}$.

The following constructions are standard. (See for example \cite{Thurston}). For  $R>0$, let 
\[    {\tilde M}_{R}^{\pm}~ =~ \{ (X,Y,Z)~ |~  X^2+Y^2\pm Z^2 = \pm R^2 \}, \] 
and $M_R^{\pm}$ the connected component of ${\tilde M}_R^{\pm}$ that contains 
the point $(0,0,R).$
The quadratic form $X^2+Y^2 \pm Z^2$ defines a Riemannian 
metric of constant curvature $\pm R^{-2}$ on  $M_{R}^{\pm}$.
For example,  $M_{1}^-$ is isometric to 
the hyperbolic plane and $M_{1}^+$ is the unit sphere. 

We construct a projective (Klein) model as follows: 
Given a point $P \in M^{\pm}$, there exists a unique
pair $(X,Y)$ such that $(X, Y, R)$ lies in the line
containing $P$ and the origin.  This defines an 
injective map $F^-_R: M^-_R \rightarrow {\mathbb R}^2$
and a two-to-one map  $F^+_R: M^+_R \rightarrow {\mathbb R}^2$.
From henceforth, we will restrict $F^+_R$ to the upper hemisphere
of $M^+_R$.  By pulling the constant curvature metrics back by 
the diffeomorphism $(F^{\pm}_R)^{-1}$ and by setting $\kappa = \pm R^{-2}$, 
we obtain the following model for $M_{\kappa}$:

\begin{prop} \label{Model}
Let $(\rho, \theta)$ be polar coordinates on ${\mathbb R}^2$
and let
\begin{equation}  \label{gkappa}
 g_{\kappa}~ =~ 
\frac{1}{(1 +\kappa \cdot \rho^2)^2}~ d\rho^2~ +~
\frac{\rho^2}{1 +\kappa \cdot \rho^2}~ d\theta^2.
\end{equation}
If $\kappa < 0$, then $g_{\kappa}$ is a complete
Riemannian metric of constant curvature $\kappa$ on the 
Euclidean disc of radius $R=|\kappa|^{-\frac{1}{2}}$. 
If $\kappa \geq 0$, then $g_{\kappa}$ is a complete
Riemannian metric of constant curvature $\kappa$ on the 
Euclidean plane. 
\end{prop}

\begin{remk}  \label{Straight} 
Note that a path is geodesic with respect to $g_{\kappa}$ 
if and only if it is a Euclidean line segment. 
(Indeed, each geodesic in $M^{\pm}_R$ is the intersection of 
a plane through the origin and $F^{\pm}_R$ is defined
via radial projection).  In particular, given 
$\kappa, \kappa' \in {\mathbb R}$, a curve in ${\mathbb R}^2$ 
is a geodesic polygon  with respect to  $g_{\kappa}$ 
if and only if it is a polygon with respect  to  $g_{\kappa'}$
and the obvious containment conditions are satisfied. 
\end{remk}

Let $P$ be a Euclidean polygon belonging to 
the Euclidean ball $B(0,R)$.
Let $\nabla_{\kappa}$ (resp.  $dv_{\kappa}$) 
denote the gradient (resp. volume form) 
associated to $g_{\kappa}$. Define the quadratic
form 
\[   q_{\kappa}(u)\,=\,\int_{P} |\nabla_{\kappa} u|^2~ dv_{\kappa}. \]
Let $\Delta_{\kappa}$ denote the 
selfadjoint operator associated with $q_{\kappa}$ 
defined either on $H^1_0(P)$ (Dirichlet boundary condition) 
or on $H^1(P)$ (Neumann boundary condition).

 \begin{prop} \label{CurvatureAnalytic} 
 The eigenvalues of $\Delta_{\kappa}$ 
 vary analytically for $\kappa \in ( -R^{-2}, \infty)$.
 \end{prop}
 
 \begin{proof} 
For each compact interval in $(- R^{-2}, \infty)$, uniform
estimates on the coefficients of $g_{\kappa}$ on $P$ 
provide $C>0$ such that  $C^{-1} |\nabla_0 f|_0 \leq  
|\nabla_{\kappa} f|_{\kappa} \leq C |\nabla_0 f|_0$ and 
$C^{-1} \leq dv_{\kappa}/ dv_{0} \leq C$.
It follows that the form domain of $q_{\kappa}$ is independent 
of $\kappa$.  Inspection of (\ref{gkappa}) shows that 
the metric $g_{\kappa}$ depends analytically on $\kappa
\in (-R^{-2}, \infty)$. It follows that both $q_{\kappa}(u)$ and the 
quadratic form $u \rightarrow \int |u|^2 dv_{\kappa}$
depend analytically on $\kappa$.
Kato-Rellich theory thus applies, yielding the proposition.
 \end{proof}

Theorem \ref{TheoremConstCurvPoly} will follow from
 
\begin{thm}
For any $\kappa,$ almost every simply connected, 
geodesic polygon in $M_{\kappa}$ 
has simple Neumann (resp. Dirichlet) spectrum. 
\end{thm}
\begin{proof}
By Remark \ref{Straight} we can identify the 
set of geodesic $n$-gons 
in $M_{\kappa}$ with a set of Euclidean polygons. 
In particular, the set of all geodesic polygons in $M_{\kappa}$
with $\kappa \leq 0$ can be identified with the set 
$U \subset {\mathcal P_n} \times (- \infty, 0]$ defined by
\[   U~ =~  \left\{ (P, \kappa)~ |~  \kappa  \leq 0, \ P 
         \subset B \left(0,|\kappa|^{-\frac{1}{2}} \right) \right\}. \] 
By Theorem \ref{TheoremSimple} and Proposition 
\ref{CurvatureAnalytic},  for almost every $P \in {\mathcal P}_n$, 
the set of $\kappa$ such that $(P, \kappa) \in U$ has simple
spectrum is at most countable.  (See the proof of Corollary 
\ref{LemmaAnalyticSimple}.)  It follows that, with respect to 
Lebesgue measure on ${\mathcal P}_n\times \R$, the spectrum of  
almost every $(P, \kappa) \in U$ is simple.   

Given $\kappa\leq 0$, let $A_{\kappa}$ be the 
set of $P \in {\mathcal P}_n$ such that $(P, \kappa) \in U$ 
does not have simple spectrum.  Then it follows from above 
that for almost every $\kappa \in {\mathbb R}$, the set 
$A_{\kappa}$ has measure zero.  

To conclude that the measure of $A_{\kappa}$ equals
zero for every $\kappa$, we use the fact that if a metric 
$g$ is rescaled by a constant $k$, then the Laplacian is rescaled
by $k^{-1}$. In particular, the multiplicity of each eigenvalue does 
not change if the metric is rescaled. 

On the other hand, the curvature of the metric does change 
under rescaling if the initial curvature is nonzero. 
In particular, the metric $(\kappa/ \kappa') \cdot  g_{\kappa}$ 
on $B(0, |\kappa|^{-\frac{1}{2}})$ has curvature $\kappa'$
and serves as another model for $M_{\kappa'}$.
It follows that if $A_{\kappa}$ has measure zero
then so does $A_{\kappa'}$.  

The case of nonnegative curvature is proven similarly.
\end{proof}



\hfill \\
\hfill \\

\noindent
 \\
Laboratoire de Math\'ematiques Jean Leray\\
UMR CNRS 6629-Universit\'e de Nantes\\
2 rue de la Houssini\`ere\\
BP 92 208\\
F-44 322 Nantes Cedex 3\\ 
{\tt Luc.Hillairet@math.univ-nantes.fr}\\

\noindent
Department of Mathematics, \\
Indiana University, Bloomington, IN, 47401 \\
{\tt cjudge@indiana.edu} 

\end{document}